
\newcount\capitolo
\newcount\incremento
\newcount\paragrafo

 3

\font\bla=msym10
\font\eightbf=cmbx8

\def\acapo{ {\vskip  0mm plus .2mm } \noindent }
 3
\def\refe#1{\par\bigbreak\centerline{\bf{\bfsmc#1}}
\par\medbreak}
\def\bff#1{{\bf #1}}
\def\bfsmc{\eightbf\bff}
\def\para#1{\global\advance\paragrafo1%
\par\bigbreak\centerline{\bf\number\paragrafo.\ \bfsmc#1}%
\global\incremento=0%
\par\medbreak}
\def\df{\global\advance\incremento by1
{\bf Definition\ \ }}

\def\pf{{\it Proof.\ \ }}
\def\prop{\global\advance\incremento by1 {\bf Proposition
\ \ }}
\def\th{\global\advance\incremento by1 {\bf Theorem
\ \ }}
\def\rem{\global\advance\incremento by1 {\bf Remark
\ \ }}
\def\lem{\global\advance\incremento by1 {\bf Lemma
\ \ }}
\def\cor{\global\advance\incremento by1 {\bf Corollary
\ \ }}

\def\vp{\varphi}

\def\q{\hbox{\bla Q}}
\def\qu{{\cal Q}}

\def\B{{\rm B}}
\def\C{{\hbox{\bla C}}}

\def\L{{\rm L}}
\def\R{{\rm R}}
\def\S{{\rm S}}

\def\V{{\rm V}}

\def\co{{\Delta}}

\def\pr{^{\prime}}

\def\Y{{\rm Y}}

\def\T{{\rm T}}
\def\a{{\alpha}}
\def\b{{\beta}}
\def\e{{\varepsilon}}
\def\ga{{\gamma}}

\def\x{{\Lambda}}
\def\w{{\wedge}}
\def\ot{{\otimes}}
\def\mo{{L_{\l/\mu}\V}}
\def\sca{{\sqcap\hskip-6.5pt\sqcup}}
\def\l{{\lambda}}
\def\o{{\omega}}

\def\s{{\cal S}}
\def\vp{{\varphi}}
\def\d{{\delta}}
\def\si{{\sigma}}
\def\vir{,\ldots ,}

\def\lra{{\longrightarrow}}

\def\P {{\rm P}}
\def\cvd{\hfill $\sqcap \hskip-6.5pt \sqcup$}  
\footline={\hss\lower .8cm \hbox{\tenrm\folio}\hss} 
\hsize=15 true cm
\vsize=20.5 true cm
\baselineskip=.5cm   

\parindent=0pt
\topskip=0cm
\voffset 1.5 truecm 
\binoppenalty=10000  
\relpenalty=10000    
\def\interi{\hbox{\bla Z}}
\def \nat{\hbox{\bla N}}

{\bf A LITTLEWOOD-RICHARDSON FILTRATION AT ROOTS OF 1 \acapo
FOR MULTIPARAMETR DEFORMATIONS OF SKEW SCHUR MODULES
\vskip 10mm
G.BOFFI$^\star$ and M.VARAGNOLO}$^\star$\footnote{}{$^\star$ Dipartimento di
Matematica, Universit\`a di Roma Tor Vergata, Via della Ricerca Scientifica,
00133 Roma, Italy. E-mail address : boffi@vax.mat.utovrm.it,
varagnolo@vax.mat.utovrm.it}

\vskip 18mm

Let $R$ be a commutative ring, $ q$  a unit of $R$
and $\P$  a multiplicatively antisymmetric
matrix with coefficients  which are integer powers of $q$.
Denote by $SE(q,\P)$ the multiparameter quantum matrix
bialgebra associated to $q$ and
P. Slightly generalizing [H-H],
we define a multiparameter deformation $L_{\l/\mu}\V_\P$ of the classical
skew Schur module.
In  case $R$ is a field and $q$ is not a root of 1, arguments
like those given in [H-H] show that  $L_{\l/\mu}\V_\P$ is
completely reducible, and its decomposition into irreducibles
is $\sum_\nu{\ga(\l/\mu;\nu)L_\nu\V_\P},$ where
 the coefficients
$\ga(\l/\mu;\nu)$ are the usual Littlewood-Richardson
coefficients.
When $R$ is {\sl any} ring and $q$ is allowed to be a root of 1,
we construct a filtration of  $L_{\l/\mu}\V_\P$ as an
$SE(q,\P)$-comodule, such that its associated graded object
is precisely $\sum_\nu{\ga(\l/\mu;\nu)L_\nu\V_\P}.$\acapo

\para{The Ingredients}
{\bf 1.1} Let $N>1$ be a positive integer.
Choose a unit $q$ in a commutative ring $R$; fix
a  matrix $\P=(p_{ij})_{i,j=1}^N$ where the $p_{ij}$'s are non-zero
elements of
$R$ with the property
$$p_{ij}p_{ji}=p_{ii}=1\ \ \forall i,j=1\vir N.$$
Consider the free $R$-module $\V_\P$ with basis $\{u_1\vir u_N\}$
and define an automorphism $\b_{q,\P}$ on $\V_\P\otimes\V_\P$ by the
following rule:
$$\b_{q,\P}(u_i\otimes u_j)=
\cases{u_i\otimes u_i&   if $i=j$\cr
qp_{ji}u_j\otimes u_i&  if  $i<j$\cr
qp_{ji}u_j\otimes u_i +(1-q^2)u_i\otimes u_j&  if $i>j$\cr}. $$
Then $(\V_\P,\b_{q,\P})$ is a YB pair in the sense of [H-H].
 Moreover it satisfies the Iwahori's quadratic equation
$$(id_{\V_\P\otimes \V_\P}-\b_{q,\P})\circ(id_{\V_\P\otimes \V_\P}+
q^{-2}\b_{q,\P})=0,$$
as we can easily verify.\acapo
\vskip 5mm
{\bf 1.2} The {\sl multiparameter quantum matrix bialgebra} $SE(q,\P)$ [S]
is the algebra generated by the $N^2$ elements $x_{ij}\ (i,j=1\vir N)$
with relations (for $i<j$ and $k<m$) :
$$x_{ik}x_{im}=qp_{mk}x_{im}x_{ik}\ \  \ \
x_{ik}x_{jk}=qp_{ij}x_{jk}x_{ik}\ \ \ \
p_{mk}x_{im}x_{jk}=p_{ij}x_{jk}x_{im}$$
$$p_{km}x_{ik}x_{jm}-p_{ij}x_{jm}x_{ik}=(q-q^{-1})x_{im}x_{jk}.$$
The coalgebra structure is given by the following comultiplication and
counity :
$$\co(x_{ij})=\sum_{k=1}^N{x_{ik}\otimes x_{kj}}\ \ \ ,\ \ \
\e(x_{ij})=\d_{ij}.$$ \vskip 5mm
{\bf 1.3} There is a natural
 $SE(q,\P)$-comodule structure on $\V_\P$ given by
 $$u_j\mapsto \sum_i{u_i\otimes x_{ij}}.$$
Consider the  ideal ${\cal B}^+_\P$ of $SE(q,\P)$ generated by all
$x_{ij}$ with  $i>j$ and put
$$SB^+(q,\P)=SE(q,\P)/{\cal B}^+_\P.$$
 The
relations between the generators in $SB^+(q,\P)$ are those given in (1.2) when
we put $x_{ij}=0$ for $i>j$. In particular  $x_{ii}$ commutes with $x_{jj}$ for
all $i,j$.
\vskip 5mm
{\bf 1.4} Henceforth the $p_{ij}$'s  will be integer
powers of $q$. More precisely (cf.[R]) we shall take
$$p_{ij}=q^{2(u_{ji}-u_{j-1i}-u_{ji-1}+u_{j-1i-1})},$$
 where
$U=(u_{ij})_{i,j=1}^{N-1}$ is an appropriate alternating integer matrix.
 In this way we shall be in the situation of
[C-V 1-2], where in fact  an integer form of the multiparameter
quantum function algebra is constructed.
 From now on, we shall skip all indices $q,\P$  in our
notations as long as no ambiguity is likely.\acapo
\vskip 5mm
{\bf 1.5} We now begin reviewing some results of [H-H], freely adopting
the notations in there.
Starting from the YB pair $(\V,\b_\V)$, we can construct some graded
YB bialgebras. First of all the tensor algebra $T\V=\bigoplus_{i\geq 0}
{\V^{\otimes i}}=\bigoplus_{i\geq 0}T_i\V$ with  YB operator
$T(\b_\V)=\bigoplus_{i,j\geq 0} {\b_\V(\chi_{ij})}$, where $\chi_{ij}$ is the
following element of $\s_{i+j}$ :
$$\chi_{ij}=\left(\matrix{1&2&\ldots&i&i+1&i+2&\ldots&i+j\cr
j+1&j+2&\ldots&j+i&1&2&\ldots&j\cr}\right).$$
We recall that, if $\si=\si_{i_1}\cdots\si_{i_r}$ is a reduced expression
for an element $\si\in\s_k$, then it is well defined on $T_k\V$
the operator $\b_\V(\si)=\b_\V(\si_{i_1})\circ\cdots\circ \b_\V(\si_{i_r})$,
$\b_\V(\si_j)$ being the map $id_\V^{\otimes j-1}\otimes \b_\V\otimes
id_\V^{\otimes k-j-1}.$ In order to describe the coproduct of $T\V$, for every
sequence $\a= (\a_1\vir\a_s)$ of nonnegative integers with $\sum_i\a_i=k$,
define $\co^\a_ {T\V}$ to be  the composite map $T\V\lra T_s\V\lra T_\a\V$ of
the
$s$-th iteration of $\co_{T\V}$ and the projection onto $T_\a\V=
\V^{\otimes \a_1}\otimes\cdots\otimes\V^{\otimes\a_s}.$
Put
$$\s^\a=\{\si\in\s_k | \si(1)<\cdot\cdot
<\si(\a_1),\si(\a_1+1)<\cdot\cdot<\si(\a_1+\a_2),\ldots,
\si(\sum_{i=1}^{s-1}{\a_i}+1)
<\cdot\cdot<\si(\sum_{i=1}^s\a_i)\}.$$
Then $\co_{T\V}^\a=\sum_{\si\in\s^\a}{\b_\V(\si^{-1})}.$\acapo
\vskip 5mm
{\bf 1.6} We consider  the symmetric and the exterior algebras $S\V$ and
$\x\V$ of the YB pair $(\V,\b_\V)$, which will play a key role in what follows.
 The
algebra $S\V$ is  generated by $u_1\vir u_N$ with relations
$$u_i u_j=p_{ji}qu_j u_i,$$
while $\x\V$ is the algebra on the same generators with relations
$$u_i\w u_i=0,\ \ \ p_{ji}qu_i\w u_j+u_j\w u_i=0\ (i<j).$$
So for every sequence $i=(i_1\vir i_k)$ of elements in [$1,N$] we have
$$u_{i_1}\w\cdots\w u_{i_k}=
\cases{0&if there are repetitions in $i$\cr
(\prod_{r<t,\si(r)>\si(t)}-q^{-1}p_{i_{\si(r)}i_{\si(t)}})u_{i_{\si(1)}}\w
\cdots\w u_{i_{\si(r)}}&if $i_1<\cdots<i_k$ and $\si\in\s_k$\cr}.$$
The $R$-modules $S_r\V$ and $\x_r\V$ are free with bases, respectively,
$$\{u_{j_1}\cdots u_{j_r}|\ 1\leq j_1\leq\cdots\leq j_r\leq N\}
\ \ ,\ \  \{u_{j_1}\w\cdots\w u_{j_r}|\ 1\leq j_1<\cdots< j_r\leq N\}.$$
\vskip 5mm
{\bf 1.7} Put $\ga_\V=-q^{-2}\b_\V.$ Then the two YB operators $\b_\V,\ga_\V$
satisfy conditions (4.9) and (4.10) in [H-H], that is, $(\V,\b_\V,\ga_\V)$
is a YB triple. From this follows (Theorem 4.10 in [H-H]) that $S\V$ and $\x\V$
are graded YB bialgebras. Moreover there exist YB operators $\varphi_
{S\V}$, $\psi_{S\V}$ on $S\V$, and  $\varphi_
{\x\V}$, $\psi_{\x\V}$ on $\x\V$, for which $(S\V, \varphi_
{S\V}, \psi_{S\V})$ and  $(\x\V, \varphi_
{\x\V}, \psi_{\x\V})$ are YB algebra triples. In particular, the operator
 $ \varphi_{\x\V}$ is defined by the relation  $\varphi_{\x\V}\circ(p\ot p)=
(p\ot p)\circ T(-\b_\V)$ where $p$ denotes the projection from $T\V$ onto
$\x\V$. The multiplicative structure on $\x\V$
is given by the fusion procedure, namely, by
$$m_{T_i(\x\V)}=m_{\x\V}^{\ot i}\circ \varphi_{\x\V}(\o_i),\ \ \ \
\o_i=\left(\matrix{1&2&\cdots&i&i+1&i+3&\cdots&2i\cr
1&3&\cdots&2i-1&2&4&\cdots&2i\cr}\right).$$
Finally note that $T\V$, $S\V$ and $\x\V$ are $SE$-equivariant
as YB bialgebras with  YB algebra triples, that is, all the structure morphisms
(including YB operators) are homomorphisms of
$SE$-comodules.\acapo \vskip 5mm
{\bf 1.8} A translation into our setting of Lemma 5.3 in [H-H] gives
the following very useful equality.\acapo
\vskip 3mm
\lem {\sl For any $k\geq 0$ and any sequence $(i_1\vir i_k)$ with
$1\leq i_1<\cdots<i_k\leq N$ we have :
$$\co^{(1\vir 1)}_{\x\V}(u_{i_1}\w\cdots\w u_{i_k})=
\sum_{\si\in\s_k}{(\prod_{r<t, \si(r)>\si(t)}-qp_{i_{\si(r)}i_{\si(t)}})
u_{i_{\si(1)}}\ot\cdots\ot u_{i_{\si(k)}}}.$$
In particular, for any $k$, $\co_{\x\V}:\x\V\lra T_k\V$ is a split injection.}
\cvd\acapo
\vskip 5mm
{\bf 1.9} We are now ready to introduce our multiparameter
deformations of Schur modules. In fact all definitions and results in
Section 6 of [H-H], stated for the "Jimbo case", still hold in our situation.
For all but  Lemma 6.12 can be deduced from formal properties of graded YB
bialgebras which are also equipped with a structure of YB algebra triple. The
proof of Lemma 6.12, which depends directly on
the definition of $\b_\V$, can be easily
 modified for our purposes.\acapo
Given a skew partition $\l/\mu$ with $l(\l/\mu)=s$ and $\l_1=t$, denote
by $d_{\l/\mu}(\V)$ the {\sl Schur map}, that is, the composite map
$$\x_{\l/\mu}\V=\x_{\l_1-\mu_1}\V\ot\cdots\ot \x_{\l_t-\mu_t}\V
\buildrel{{\scriptstyle
\co_{\x\V}^{(1^{\l_1-\mu_1})}\ot\cdots
\ot \co_{\x\V}^{(1^{\l_t-\mu_t})}}}\over\lra T_{\l/\mu}\V=
T_{\l_1-\mu_1}\V\ot\cdots\ot T_{\l_t-\mu_t}\V\lra$$
$$\buildrel{{\scriptstyle (-q^{-2 }\b_\V)(\chi_{\l/\mu})}}\over\lra
T_{\tilde{\l}/\tilde{\mu}}\V=T_{\tilde{\l_1}-\tilde{\mu_1}}\V\ot\cdots
\ot T_{\tilde{\l_s}-\tilde{\mu_s}}\V\buildrel{{\scriptstyle p\ot\cdots\ot p}}
\over\lra
S_{\tilde{\l}/\tilde{\mu}}\V=S_{\tilde{\l_1}-\tilde{\mu_1}}\V\ot\cdots
\ot S_{\tilde{\l_s}-\tilde{\mu_s}}\V,$$
where, as usual, $\tilde{\l}$ denotes the dual partition of $\l$, and $\chi_
{\l/\mu}$ is the permutation defined in Section 6 of [H-H].
We illustrate such a permutation by the following example :
 $$\l=(5,4,2)\ \ \ \ \mu=(2,1)\ \
\ \ \chi_{\l/\mu}= \left(\matrix{1&2&3&4&5&6&7&8\cr4&6&8&2&5&7&1&3\cr}\right)$$
$$\matrix{\bullet&\bullet&1&2&3\cr
\bullet&4&5&6& \cr
7&8& & & \cr}\ \ \
\buildrel{\chi_{\l/\mu}}\over\lra\ \ \
\matrix{\bullet&\bullet&4&6&8\cr
\bullet&2&5&7& \cr
1&3& & & \cr}.$$
The image of the Schur map, denoted by $\mo$, is the
{\sl Schur module of $\V$ with respect to the skew partition $\l/\mu.$}
It is an $SE$-comodule, with coaction induced by the following coaction
 on $T_k\V$ :
$$u_{j_1}\ot\cdots\ot u_{j_k}\mapsto\sum_{i_1\vir i_k}
{(u_{i_1}\ot\cdots\ot u_{i_k})\ot x_{i_1j_1}\ot\cdots\ot x_{i_kj_k}}.$$
\vskip 5mm
{\bf 1.10} The principal properties of $\mo$ are summarized  in the following
theorem, which one proves along the lines of Theorem 6.19 and Corollary
6.20 in [H-H].\acapo \vskip 3mm
\th {\sl Let $\l/\mu$ a skew partition with $l(\l/\mu)=s$.
Then :\acapo (i) $\mo$ is an $R$-free module, and for any
$\si\in\s_N$, a free basis is the set
 $$L_{\l/\mu}\Y(\si)=\{d_{\l/\mu}(\V)(\xi_\S)|\ \S\in
St_{\l/\mu}\Y(\si)\}.$$
Here $St_{\l/\mu}\Y$ denotes the set of all standard tableaux in the
alphabet $\Y(\si)=\{u_{\si(1)}<\cdots<u_{\si(N)}\},$ and
$$\xi_\S=
\S(1,\mu_1+1)\w\cdots\w\S(1,\l_1)\ot\cdots\ot\S(s,\mu_s+1)\w\cdots\w
\S(s,\l_s)\in\x_{\l/\mu}\V.$$
(ii) Let $R\pr$ be a commutative ring and let  $f:R\lra R\pr$ be a
homomorphism of commutative rings. Then we have an isomorphism of
$SE\pr$-comodules
$$L_{\l/\mu}(R\pr\ot_R\V)\simeq R\pr\ot_RL_{\l/\mu}\V\ \ ,
\ \ E\pr=R\pr\ot_RE.$$}
\cvd
\vskip 3mm
As a consequence of $(ii)$, it will not be  restrictive for us to take
$R=\interi[\qu,\qu^{-1}]$, where $\qu$ stands for an indeterminate.\acapo
 \vskip 5mm
{\bf 1.11} We recall that an element  of
$Tab_{\l/\mu}\Y(\si)$, the set of all tableaux of shape $\l/\mu$ with elements
in $\Y(\si)$, is said to be {\sl row-standard} if its rows are strictly
increasing, and {\sl column-standard} if its columns are non-decreasing.
A tableau  is said to be {\sl standard} if it is both row- and column-standard.
 Let $Row_{\l/\mu}\Y(\si)$ denote the set of
row-standard tableaux  of shape $\l/\mu$  and with elements
in  $\Y(\si)$. For every $\S\in Row_{\l/\mu}\Y(\si)$, the element
$d_{\l/\mu}(\V)(\xi_\S)$ can be expressed as a linear combination of basis
elements. The algorithm, call it ${\cal R}_\si$, which does this is based on a
descending induction with respect to a pseudo order defined in
$Tab_{\l/\mu}\Y(\si)$.
 Let $\S$ and $\S\pr$ be elements
in $Tab_{\l/\mu}\Y(\si).$
 We say that $\S\leq_\si\S\pr$ if
$\forall p,q$\acapo
 $\#\{(i,j)\in\co_{\l/\mu}|\ i\leq p,\ \S(i,j)\in\{u_{\si(1)}\vir
u_{\si(q)}\}\}\geq$\acapo
\hskip 70mm$ \#\{(i,j)\in\co_{\l/\mu}|\ i\leq p,\
\S\pr(i,j)\in\{u_{\si(1)}\vir u_{\si(q)}\}\}.$
The key steps of ${\cal R}_\si$
are the following:\acapo
1. Choose two adjacent lines in S where there is a violation of
column-standardness; we are in the situation of Proposition
(1.12) below, and we can use Corollary (1.13). We get certain $\S_i$'s
such that $\S_i<\S$ for every $i$.\acapo
2. Reorder in increasing order
$\S_i(1,\mu_1+1)\w\cdots\w\S_i(1,\l_1)\vir
\S_i(s,\mu_s+1)\w\cdots\w\S_i(s,\l_s)$ for each $i$; this operation produces a
power of $q$
 for every $\S_i$ (cf. (1.6)).\acapo
3. Apply induction to each $\S_i$.\acapo
${\cal R}_\si$ is also called the "straightening law with respect to the
ordering $u_{\si(1)}<\cdots<u_{\si(N)}$".
\vskip 5mm
{\bf 1.12} \prop {\sl Let $\l=(\l_1,\l_2)$ and $\mu=(\mu_1,
\mu_2)$ be partitions with $\l\supset\mu$. Define $\ga=\l-\mu$ and take
$a,b$ nonnegative integers with $a+b<\l_2-\mu_1$. Then the image of the
composite map
$$\overline\sca_{(a,b)}:\x_a\V\ot\x_{\ga_1-a+\ga_2-b}\V\ot\x_b\V
\buildrel{1\ot\co\ot 1}\over\lra
\x_a\V\ot\x_{\ga_1-a}\V\ot\x_{\ga_2-b}\V\ot\x_b\V
\buildrel{m\ot m}\over\lra
\x_{\ga_1}\V\ot\x_{\ga_2}\V=\x_{\l/\mu}\V
$$
is contained in $Im(\sca_{\l/\mu})$, where $\sca_{\l/\mu}$ is given by
$$\sum_{\nu=0}^{\l_2-\mu_1}
{\x_{\ga_1+\ga_2-\nu}\V\ot\x_\nu\V}
\buildrel{\co\ot 1}\over\lra
\sum_{\nu=0}^{\l_2-\mu_1}{\x_{\ga_1}\V\ot\x_{\ga_2-\nu}\V\ot\x_\nu\V}
\buildrel{1\ot m}\over\lra
\sum_{\nu=o}^{\l_2-\mu_1}{\x_{\ga_1}\V\ot\x_{\ga_2}\V}.$$}
\vskip 3mm
\pf Mimic the proof of Lemma 6.15 in [H-H].\cvd
\vskip 5mm
{\bf 1.13} \cor {\sl Let $\l/\mu$ be a skew partition with
$l(\l)=s$,  $\si$ be an element of $\s_N$ and  S be an element of $
Row_{\l/\mu}\Y(\si)\setminus St_{\l/\mu}\Y(\si)$. Then there exist $\S_1\vir
\S_r\in Row_{\l/\mu}\Y(\si)$
 ($r\in\nat $) with $\S_i<_{\si}\S,\ \forall i=1\vir r$ such that
$$\xi_\S-\sum_ic_i\xi_{\S_i}\in Im(\sca_{\l/\mu})=
Ker (d_{\l/\mu}(\V)),$$
for some $c_i\in\interi[q,q^{-1}]$. Here :
$$\sca_{\l/\mu}=
\sum_{i=1}^{s-1}{1_1\ot\cdots\ot 1_{i-1}\ot\sca_{\l^i/\mu^i}\ot
1_{i+2}\ot\cdots
\ot 1_s},\  \l^i=(\l_i,\l_{i+1}),\ \mu^i=(\mu_i,\mu_{i+1}),\
1_j=id_{\x_{\l_j-\mu_j}\V} .$$}
\vskip 3mm
\pf Mimic the proof of Lemma 6.18 in [H-H].\cvd
\vskip 5mm
{\bf 1.14} We want to stress a consequence of Theorem (1.10)
 and of all the machinery which allows to prove it. First of all
note that the subcategory of ${\cal YB}_R$ (cf. [H-H]) given by the
YB pairs as
in  1.1 is a preadditive one. Namely, let $\P^1=(p^1_{ij})_{i,j=1}^n$ and
$\P^2=(p^2_{ij})_{i,j=1}^m$ be two multiplicatively antisymmetric matrices,
 and
put $\V_{\P^1}=<u_1^1\vir u_n^1>$, $\V_{\P^2}=<u_1^2\vir u_m^2>.$
We define a YB operator on $\V_{\P^1}\oplus\V_{\P^2}$ by means of the
matrix $\P=(p_{ij})_{i,j=1}^N$, $N=n+m$, defined as follows:
$$p_{ij}=\cases{p^1_{ij}& for $i,j\in [1,N]$\cr
p_{ij}^2& for $i,j\in [n+1,N]$\cr
1& for $i\in[1,n],j\in [n+1,N]$ or $i\in [n+1,N],j\in [1,n]$\cr}.$$
Then $\b_\P$ is a YB operator on $\V_\P=<u_1^1\vir u_n^1,u_1^2\vir u_m^2>.$
Note that $\V_\P$ becomes in a natural way an $SE(q,\P^1)\ot
SE(q,\P^2)$-comodule.\acapo
 Write for short $\V_i=\V_{\P^i}$, $\b_i=\b_{\P^i}$,
 for
$i=1,2$, and let $\mu \subset\ga\subset\l$ be partitions.
Following [A-B-W],
define two  R-modules
$$M_\ga(\x_{\l/\mu}(\V_1\oplus\V_2))=Im(
\sum_{\mu\subseteq\si\subseteq\l,\si\geq \ga}
{\x_{\si/\mu}\V_1\ot\x_{\l/\si}\V_2\lra\x_{\l/\mu}(\V_1\oplus\V_2))},$$
$$\dot{ M}_\ga(\x_{\l/\mu}(\V_1\oplus\V_2))=Im(
\sum_{\mu\subseteq\si\subseteq\l,\si> \ga}
{\x_{\si/\mu}\V_1\ot\x_{\l/\si}\V_2\lra\x_{\l/\mu}(\V_1\oplus\V_2))},$$
where the indicated maps are obtained by tensoring the obvious maps
$$\x_{\si_i-\mu_i}\V_1\ot\x_{\l_i-\si_i}\V_2\lra\x_{\l_i-\mu_i}(\V_1\oplus
\V_2).$$
Let $M_\ga(L_{\l/\mu}(\V_1\oplus\V_2))$ and
$\dot{M}_\ga(L_{\l/\mu}(\V_1\oplus\V_2))$ be the images of the previuos modules
under the Schur map $d_{\l/\mu}(\V_1\oplus\V_2)$.
The following result holds as in the classical case :\vskip 3mm
\th  {\sl
The $R$-modules
$$L_{\l/\mu}\V_1\otimes L_{\l/\ga}\V_2\ \ ,\ \
M_\ga(L_{\l/\mu}(\V_1\oplus\V_2))/\dot{M}_\ga(L_{\l/\mu}(\V_1\oplus\V_2))$$
are isomorphic. Hence the $R$-modules $M_\ga(L_{\l/\mu}(\V_1\oplus\V_2))$,
$\mu\subseteq\ga\subseteq\l$,
give a filtration of $L_{\l/\mu}(\V_1\oplus \V_2)$,
whose associated graded module is isomorphic to
$$\sum_{\mu\subseteq\ga\subseteq\l}L_{\ga/\mu}\V_1\otimes L_{\l/\ga}\V_2.$$}
\vskip 3mm \pf Follow {\sl verbatim} the proof of Theorem II. 4.11 in
[A-B-W].\cvd
\vskip 3mm
Note that the isomorphism of the theorem is in fact an isomorphism of
$SE(q,\P^1)\ot SE(q,\P^2)$-comodules.\acapo
\vskip 5mm
\vfill\eject
\para{ The Recipe}
{\bf 2.1} In this Section we let $R$ be the ring $R=\interi[\qu,\qu^{-1}]$,
$\qu$ an indeterminate, and take a multiplicatively antisymmetric matrix
$\P=(p_{ij})_{i,j=1}^N$, and the YB pair $(\V_\P,\b_{\qu,\P})$, where
$\V_\P=<u_1\vir u_N>$ and
$$\b_{\qu,\P}(u_i\ot u_j)=
\cases{u_i\ot u_i& if $i=j$\cr
\qu p_{ji}u_j\ot u_i&if $i<j$\cr
\qu p_{ji}u_j\ot u_i+(1-\qu^2)u_i\ot u_j& if $i>j$\cr}.\leqno(1)$$
We are going to construct a filtration of $L_{\l/\mu}\V_\P$ as an
$SE(\qu,\P)$-comodule, such that the associated graded object
is isomorphic to $\sum_\nu{\ga(\l/\mu;\nu)L_\nu\V_\P}.$ As in the
classical Littlewood-Richardson rule, here $\ga(\l/\mu;\nu)$
stands for the number of standard tableaux of shape $\l/\mu$
filled with $\tilde{\mu}_1$ copies of $1$, $\tilde{\mu}_2$ copies
of $2$, $\tilde{\mu}_3$ copies of $3$, etc., such that the
associated word (formed by listing all entries from bottom to
top in each column, starting from the leftmost column) is a
lattice permutation.\acapo
The construction is  a suitable "deformation" of the one used
in the first author's doctoral thesis, Brandeis University 1984,
as illustrated for instance in [B]. We again remark that owing
to Theorem (1.10) (ii), the construction holds in fact for every
commutative
ring $R$ and every choice of a unit $q\in R$.\acapo
\vskip 5mm
{\bf 2.2} In order to embed $L_{\l/\mu}\V_\P$ into a (non-skew)
Schur module, let $M=\mu_1$ and consider another multiplicatively
antisymmetric matrix $\P\pr=(p_{ij}\pr)_{i,j=1}^M$, together with
the YB pair $(\V_{\P\pr},\b_{\qu,\P\pr})$, where $\V_{\P\pr}=
<u\pr_1\vir u_M\pr>$ and $\b_{\qu,\P\pr}$ is defined similarly to
(1) above. For convenience of notations, we shall  denote $\V_\P,
u_i,\V_{\P\pr}$, and $u_i\pr$ by $\V,i,\V\pr$, and $i\pr$,
respectively.\acapo
It follows from Theorem (1.14) that the $SE(\qu,\P\pr)
\ot SE(\qu,\P)$-comodule $L_\l(\V\pr\oplus \V)$  is isomorphic
to $\sum_{\a\subseteq\l}{L_\a\V\pr\ot L_{\l/\a}\V},$ up to a
filtration.\acapo
Let $(L_\l(\V\pr\oplus \V))_h$ denote the sub-$R$-module of
$L_\l(\V\pr\oplus \V)$ spanned by the tableaux in which $h$
$\V\pr$-indices occur. (In this section we identify tableaux
and corresponding elements of Schur modules.) Then up to a
filtration,
$$(L_\l(\V\pr\oplus \V))_h\simeq\sum_{\a\subseteq\l,|\a|=h}
{L_\a\V\pr\ot L_{\l/\a}\V},$$
as $SE(\qu,\P\pr)\ot SE(\qu,\P)$-comodules.\acapo
If $(L_\l(\V\pr\oplus \V))_{\tilde{\mu}}$ denotes the
sub-$R$-module of $L_\l(\V\pr\oplus \V)$ spanned by the
tableaux in which every $i\pr$ occurs exactly $\tilde{\mu}_i$
times, also :
$$(L_\l(\V\pr\oplus \V))_{\tilde{\mu}}\simeq
\sum_{\a\subseteq\l}{(L_\a\V\pr)_{\tilde{\mu}}\ot
L_{\l/\a}\V},\leqno(2)$$
as $SE(\qu,\P)$-comodules, up to a filtration.\acapo
Since the bottom piece of the filtration relative to (2)
corresponds to the (lexicographically) largest partition
$\a$, namely $\mu$, it follows :
$$(L_\l\V\pr)_{\tilde{\mu}}\ot L_{\l/\mu}\V
\buildrel {SE(\qu,\P)}\over\hookrightarrow
(L_\l(\V\pr\oplus\V))_{\tilde{\mu}}.$$
And $rk(L_\mu\V\pr)_{\tilde{\mu}}=1$
implies that
$$L_{\l/\mu}\V\buildrel{SE(\qu,\P)}\over
\hookrightarrow (L_\l(\V\pr\oplus\V))_{\tilde{\mu}},$$
 as
wished.\acapo
Explicitly, the embedding sends the tableau $d_{\l/\mu}(\V)
(a_1\ot\cdots\ot a_s),$ $s=l(\l),$ to
$$d_\l(\V\pr\oplus\V)[(b^{(\mu_1)}\w a_1)\ot\cdots\ot (b^{(\mu_r)}\w
a_r)\ot a_{r+1}\ot\cdots\ot a_s],\ \ r=l(\mu),$$
where we write $b^{(k)}$ for $1\pr\w 2\pr\w
\cdots \w k\pr\in\Lambda_k\V\pr.$
Notice that $b^{(k)}$ is a relative
$SB^+(\qu,\P\pr)$-invariant.\acapo
\vskip 5mm
{\bf 2.3} Let {\bf t}$=(t_{r1}\vir t_{11};t_{r2}
\vir t_{12};\ldots ;t_{rs}\vir t_{1s})$ be a family of
nonnegative integers such that
$$\sum_{i=1}^s t_{ji}=\mu_j\ \ \ \ \forall j=1\vir r.$$
Let $f$ denote the $SE(\qu,\P\pr)$-equivariant composite
map :
$$\matrix{
\Lambda_{\mu_r}\V\pr\ot \cdots \ot
\Lambda_{\mu_1}\V\pr\cr
\cr
\downarrow  {\scriptstyle  \ot_{j=r}^1(\co^{t_j}_{\Lambda\V\pr})}\cr
\cr
(\Lambda_{t_{r1}}\V\pr\ot\cdots\ot\Lambda_{t_{rs}}\V\pr) \ot
\cdots\ot(\Lambda_{t_{11}}\V\pr\ot\cdots\ot\Lambda_{t_{1s}}\V\pr)\cr
\cr
\downarrow {\scriptstyle  \vp_{\Lambda\V\pr}(\o_{rs})}\cr
\cr
(\Lambda_{t_{r1}}\V\pr\ot\Lambda_{t_{r-1,1}}\V\pr\ot\cdots
\ot\Lambda_{t_{11}}\V\pr)\ot\cdots\ot(\Lambda_{t_{rs}}\V\pr
\ot\Lambda_{t_{r-1,s}}\V\pr\ot\cdots\ot\Lambda_{t_{1s}}\V\pr)\cr
\cr
\downarrow { \scriptstyle (m_{\Lambda\V\pr}^{(r)})^{\ot s}}\cr
\cr
\Lambda_{t_{r1}+t_{r-1,1}+\cdots+t_{11}}\V\pr\ot\cdots\ot
\Lambda_{t_{rs}+t_{r-1,s}+\cdots+t_{1s}}\V\pr\cr}$$
where
$t_j=(t_{j1}\vir t_{js})$, $m_{\Lambda\V\pr}^{(r)} :\Lambda\V\pr\ot
\cdots\ot\Lambda\V\pr\lra \Lambda\V\pr$ is obtained by iterating
the multiplication, and
$$\o_{rs}={\left(\matrix{
1&2&3&\cdots&s&s+1&s+2&\cdots&2s+1&\cdots&rs\cr
1&r+1&2r+1&\cdots&(s-1)r+1&2&r+2&\cdots&3&\cdots&rs\cr}\right)}$$
(cf. items (1.5) and (1.7)).\acapo
As $b^{(\mu_r)}\ot\cdots\ot b^{(\mu_1)}$ is a relative
$SB^+(\qu,\P\pr)$-invariant, also $f(b^{(\mu_r)}\ot\cdots
\ot b^{(\mu_1)})$ is so. We denote the latter by
$b({\bf t}).$\acapo
\vskip 5mm
{\bf 2.4} For every $\nu\subseteq\l$ such that $|\nu|=|\l|-|\mu|$,
let $B(\l/\nu)$ denote the set of all possible $b({\bf t})$
which satisfy the further equalities :
$$\sum_{j=1}^rt_{ji}=\l_i-\nu_i\ \ \ \ \forall i=1\vir s.$$
For every $b\in B(\l/\nu)$, we call $\vp(\nu,b)$ the restriction to
$\Lambda_\nu\V\ot\{b\}$ of the following composite map
$$\Lambda_\nu\V\ot\Lambda_{\l/\nu}\V\pr\buildrel{\vp_\nu(\l)}
\over\lra\Lambda_\l(\V\pr\oplus\V)\buildrel{d_\l(\V\pr\oplus\V)}
\over\lra\Lambda_\l(\V\pr\oplus\V),$$
where $\vp_\nu(\l)$ is obtained by tensoring the morphisms
$$\Lambda_{\nu_i}\V\ot\Lambda_{\l_i-\nu_i}\V\pr\lra
\Lambda_{\l_i}(\V\pr\oplus\V),\ \ \ \ x\ot y\mapsto x\w y,\ \ \ \
i=1\vir s.$$
\vskip 3mm
\prop {\sl The image of $\vp(\nu,b)$ lies in
$\L_{\l/\mu}\V\hookrightarrow \L_\l(\V\pr\oplus\V).$}\acapo
\vskip 3mm
\pf As $\vp(\nu,b)$ is $SE(\qu,\P\pr)\ot SE(\qu,\P)$-equivariant,
and $b$ is a relative $SB^+(\qu,\P\pr)$-invariant of $\V\pr$-content
$\tilde{\mu}$ (i.e., it contains $\tilde{\mu}_i$ copies of $i\pr$), each
element of $Im(\vp(\nu,b))$ is a relative $SB^+(\qu,\P\pr)$-invariant
of $\V\pr$-content $\tilde{\mu}$. But then we are through, thanks to
Lemma (2.5) below and to the fact that $d_\mu(\V\pr)(b^{(\mu_1)}\ot
\cdots\ot b^{(\mu_r)})$ is the only canonical tableau of content
$\tilde{\mu}$.\cvd
\vskip 5mm
{\bf 2.5} \lem {\sl  For every partition $\a$,
take in $L_\a\V\pr\ot_\R\q(\qu)$ the element
$$C_\a=d_\a(\V\pr)(1\pr\w\cdots\w\a_1\pr\ot1\pr\w\cdots\w \a_2\pr\ot
\cdots\ot 1\pr\w\cdots\w\a_l\pr),\ \ \ \ \ l=l(\a)$$
($C_\a$ is sometimes called the "canonical tableau of $L_\a\V\pr$").
Then the relative $SB^+(\qu,\P\pr)$-invariant
elements of $L_\a\V\pr\ot_R\q(\qu)$ are spanned (over $\q(\qu)$) by $C_\a$.}
\vskip 3mm
\pf Combine $(L_\a\V\pr)_{\tilde{a}}=R\cdot  C_\a$
with a multiparameter version of a suitable analogue of Theorem
6.5.2 in [P-W].\cvd
\vskip 5mm
{\bf 2.6} For each $\nu\subseteq\l$ such that $\ga(\l/\mu;\nu)
\not=0$, we wish to describe a subset of $B(\l/\nu)$, say $B\pr(\l/\nu)$,
such that $\#\B\pr(\l/\nu)=\ga(\l/\mu;\nu).$
Let $\T\in L_{\l/\nu}\V\pr$ be a standard tableau, of content
$\tilde{\mu}$, and such that its associated word, $as(\T)=
(a_1\vir a_{|\mu|}),$ is a lattice permutation.
Then $\mu$ is the content of the transpose lattice permutation
$(as(\T))\tilde{} .$
(Explicitely, $(as(\T))\tilde{} =(\tilde{a}_1\vir \tilde{a}_{|\mu|}),$
where $\tilde{a}_i$ is the number
of times $a_i$ occurs in $as(\T)$ in the range $(a_1\vir a_i).)$
Let $\tilde{\T}$ be the tableau obtained from T by replacing every entry
$a_i$ of T by $\tilde{a}_i$. For each $i\in\{1\vir s\}$ and each $j\in
\{1\vir r\}$, we set :
$$t_{ji}=\# {\rm \  of\ } j{\rm 's\  occuring\  in\  the\  }
i-{\rm th\ row\ of\ }
\tilde{\T}.$$
We denote by $b(\T)$ the element $b({\bf t})\in
B(\l/\nu)$, corresponding to this choice of $t_{ji}$'s.\acapo
\vskip 5mm
{\bf 2.7} Given any row-standard tableau T, we can consider the word
$w(\T)$ formed by writing one after the other all the rows of T,
starting from the top. As all such words can be ordered lexicographically,
we can say that $\T<_{lex}\T\pr$ if and only if $w(\T)<_{lex}w(\T\pr).$
It is then easy to see that the following holds.\acapo
\vskip 3mm
\prop {\sl If we write $b(\T)\in\Lambda_{\l/\nu}\V\pr$ as a linear
combination of row-standard tableaux, then
$$b(\T)=\pm\qu^*\T+\sum_k{c_k\T_k},\ \ \ \ c_k\in \interi[\qu,\qu^{-1}],$$
where $\qu^*$ stands for a power of $\qu$, and each $\T_k$ is a row-standard
tableau $<_{lex}\T.$}\cvd
\vskip 3mm
Since there are exactly $\ga(\l/\mu;\nu)$ tableaux $\T\in L_{\l/\nu}\V\pr$
which are standard, of content $\tilde{\mu}$, and such that $as(\T)$ is
a lattice permutation, the above Proposition implies that the elements
$b(\T)$ form a subset of $B(\l/\nu)$ of cardinality $\ga(\l/\mu;\nu).$
It is precisely this subset which we call $B\pr(\l/\nu).$\acapo
\vskip 5mm
{\bf 2.8} Consider the family of elements of
$L_{\l/\mu}\V\hookrightarrow L_\l(\V\pr\oplus\V)$ :
$${\cal F}=\{\vp(\nu,b)(x)|\ \ga(\l/\mu;\nu)\not=0,\ b\in B\pr(\l/\nu),
\ {\rm and }\ d_\nu(\V)(x)\ {\rm is\ a \ standard\ tableau}\}.$$
We claim that ${\cal F}$ is an $R$-basis of $L_{\l/\mu}\V.$\acapo
\vskip 3mm
\prop {\sl The elements of ${\cal F}$ are linearly indipendent over $R$.}
\acapo
\vskip 3mm
\pf Suppose that there exist nonzero coefficients $r_{\nu,b,x}\in R$
such that $\sum_{\cal F}{r_{\nu,b,x}\vp(\nu,b)(x)}=0$, i.e., such that
$\sum_{\nu,b,x}{r_{\nu,b,x}d_\l(\V\pr\oplus\V)(\vp_\nu(\l)(x\ot b))}=0$
in $L_\l(\V\pr\oplus\V).$ This is the same as
$$\sum_{\nu,b}{d_\l(\V\pr\oplus\V)
(\vp_\nu(\l)(y_{\nu,b}\ot b))}=0,\leqno(3)$$
where $y_{\nu,b}=\sum_x{r_{\nu,b,x}x}.$ Let $\nu_0$ be the
(lexicographically) smallest $\nu$ occuring in (3). Order the set
$B\pr(\l/\nu_0)=\{b(\T_1)\vir b(\T_p)\}$ as follows :
$$b(\T_i)<b(\T_j)\  {\rm if\ and\ only\ if\ }
 w(\T_i)<_{lex}w(\T_j).$$
 Let $b_0=b(\T_0)$
be the highest $b(\T_i)\in B\pr(\l/\nu_0)$ occuring in
$\sum_{\nu,b}{d_\l(\V\pr\oplus\V)(\vp_\nu(\l)(y_{\nu,b}\ot b))}.$
Clearly, $d_\l(\V\pr\oplus\V)(\vp_{\nu_0}(\l)(y_{\nu_0,b_0}\ot b_0))$
is not in general a linear combination of standard tableaux of
$L_\l(\V\pr\oplus\V)$, with respect to the order $1<\cdots<N<1\pr<\cdots<M\pr,$
since violations of column-standardness may occur in $b_0$. Apply therefore
to $d_\l(\V\pr\oplus\V)(\vp_{\nu_0}(\l)(y_{\nu_0,b_0}\ot b_0))$
the straightening law of $L_\l(\V\pr\oplus\V)$ with respect to
$1<\cdots<N<1\pr<\cdots <M\pr.$ One gets (recall Proposition
(2.7)) :\acapo
\vskip 3mm
$\pm\qu^* d_\l(\V\pr\oplus\V)(\vp_{\nu_0}(\l)(y_{\nu_0,b_0}\ot \T_0))+$
(a linear combination of standard tableaux with V-shape $>\nu_0$)
$+$(a linear combination
of standard tableaux with V-shape$\ =\nu_0$ and
$\V\pr$-part $<_{lex}\T_0$ ).\acapo
\vskip 3mm
Because of our choice of $\nu_0$ and $b_0$,
(3) then implies that $d_\l(\V\pr\oplus\V)(\vp_{\nu_0}(\l)(y_{\nu_0,b_0}
\ot \T_0))=0$, i.e.,
$$\sum_x{r_{\nu_0,b_0,x}d_\l(\V\pr\oplus\V)
(\vp_{\nu_0}(\l)(x\ot \T_0))}=0.$$
But this is a linear combination of standard tableaux in $L_\l(\V\pr
\oplus\V),$ with respect to the order \acapo
$1<\cdots<N<1\pr<\cdots<M\pr,$ so that
$r_{\nu_0,b_0,x}=0$ for each $x$, which contradicts our assumption
on the coefficients $r_{\nu,b,x}.$\cvd
\vskip 5mm
{\bf 2.9} \cor {
\sl ${\cal F}$ is a basis for $L_{\l/\mu}\V\ot_R\q(\qu).$}
\vskip 3mm
\pf By definition of ${\cal F}$, $\#{\cal F}=rk( L_{\l/\mu}\V).$
By Theorem (1.10)(ii), the latter rank is constant on all rings.
So proposition (2.8) says that ${\cal F}$ is a basis for the vector space
$L_{\l/\mu}\V\ot_R\q(\qu).$\cvd
\vskip 5mm
{\bf 2.10} \cor $ {\cal F}$ {\sl is a basis for $L_{\l/\mu}\V.$}
\vskip 3mm
\pf It suffices to show that ${\cal F}$ is a system of generators
for $L_{\l/\mu}\V.$ Let $y\in L_\l(\V\pr\oplus\V)$ be any tableau of type
$$\matrix{1\pr&\cdot&\cdot&\cdot&\mu_1\pr&\circ&\circ&\circ\cr
1\pr&\cdot&\cdot&\mu_2\pr&\circ&\circ&\circ&\cr
1\pr&\cdot&\mu_3\pr&\circ&\circ&&&\cr
\cdot&\cdot&\circ&\circ&&&&\cr
\cdot&\cdot&\circ&\circ&&&&\cr
\circ&\circ&\circ&&&&&\cr
\circ&\circ&&&&&&\cr},$$
where the little circles stand for basis elements of V.\acapo
Since $y\in L_{\l/\mu}\V,$ Corollary (2.9) says that in the quotient field of
$R$, there exist (unique) coefficients $q_{\nu,b,x}$, such that
$$y=\sum_{\cal F}{q_{\nu,b,x}\vp(\nu,b)(x)}.\leqno(4)$$
To both sides of (4), apply the straightening law with respect to
$1<\cdots<N<1\pr<\cdots<M\pr.$ In the left-hand  side,
only coefficients in $R$ occur. In the right-hand side, if
$\nu_0$ denotes the smallest V-shape coupled with a nonzero
$\sum_x{q_{\nu,b,x}x},$ and $b_0=b(\T_0)$ denotes the highest
element of $B\pr(\l/\nu_0)$ (cf. ordering in the proof of
Proposition (2.8)) occurring with a nonzero $\sum_x{q_{\nu_0,b,x}x},$
we find that the term $\pm\qu^* d_\l(\V\pr\oplus\V)(\vp_{\nu_0}(\l)
(\sum_x{q_{\nu_0,b_0,x}x\ot\T_0}))$ must cancel with something in
 the left-hand side;
since each $d_{\nu_0}(\V)(x)\in L_\nu\V$ is standard, it follows that
$q_{\nu_0,b_0,x}\in R$ for every $x$.\acapo
Write next (4) as :
$$y-\sum_x{q_{\nu_0,b_0,x}\vp(\nu_0,b_0)(x)}=\sum_{(\nu,b)\not=
(\nu_0,b_0)}{\vp(\nu,b)(\sum_x{q_{\nu,b,x}x})}.\leqno(4\pr)$$
Reasoning for ($4\pr$) as done for (4), it follows that
$q_{\nu_1,b_1,x}\in R$, where $(\nu_1,b_1)$ is the pair $(\nu,b)$ coming
immediately before $(\nu_0,b_0)$ in the total ordering :\acapo
$(\nu,b)<(\nu\pr,b\pr)$ if and only if either $\nu>\nu\pr$,
 or
$\nu=\nu\pr$ and $b<b\pr$ in the ordering of $B\pr(\l/\nu)$ given in
the proof of Proposition (2.8).\acapo
Repeating the argument as many times as necessary, the proof is completed.\cvd
\vskip 5mm
{\bf 2.11} \th {\sl Up to a filtration, $L_{\l/\mu}\V\simeq
\sum_\nu{\ga(\l/\mu;\nu)L_\nu\V}$ as $SE(\qu,\P)$-comodules.}
\vskip 3mm
\pf For every $\nu$ such that $\ga(\l/\mu;\nu)\not=0$, let
$M_\nu$ denote the $R$-span (in $L_\l(\V\pr\oplus\V)$)
of all elements $\vp(\tau,b)(x)$ of ${\cal F}$ such that
$\tau\geq\nu$. Also let $\dot{M}_\nu$ denote the $R$-span of all
$\vp(\tau,b)(x)$ such that $\tau>\nu$. Clearly, we have
the isomorphism of free $R$-modules :
$$M_\nu/\dot{M}_\nu\buildrel{\psi_\nu}\over\lra
L_\nu\V\oplus\cdots\oplus L_\nu\V
\ \ \ \ (\ga(\l/\mu;\nu)\  {\rm summands}).$$
$\{M_\nu\}$ will be the required filtration, if we show that each
$\psi_\nu$ is an $SE(\qu,\P)$-isomorphism. In order to  do so,
it suffices to prove that for every fixed $b_0\in B\pr(\l/\nu)$,
and for every basis element $y\in\Lambda_\nu\V$,
$\vp(\nu,b_0)(y)-\vp(\nu,b_0)(\sum r_ix_i)\in\dot{M_\nu}$,
where $\sum r_id_\nu(\V)(x_i)$ is obtained by application
to the tableau $d_\nu(\V)(y)$ of the straightening law of $L_\nu\V.$
Notice however that $\vp(\nu,b_0)(y)\in L_{\l/\mu}\V\subseteq
L_\l(\V\pr\oplus\V)$ can be written in two ways :
$$\vp(\nu,b_0)(y)=\sum_{\cal F}{r_{\tau,b,x}\vp(\tau,b)(x)},
\leqno(5)$$
by Corollary (2.10), and
$$\vp(\nu,b_0)(y)=\sum{r_i\vp(\nu,b_0)(x_i)}+ L.C.,\leqno(6)$$
where $L.C.$ denotes a linear combination of tableaux,
standard with respect to  $1<\cdots<N<1\pr<$\acapo
$\cdots<M\pr,$
and with V-part $>\nu$. This last equality is obtained by eliminating in the
V-part of $\vp(\nu,b_0)(y)$ all violations of standardness, with respect to
$1<\cdots<N<1\pr<\cdots<M\pr.$\acapo
Comparing (5) and (6), it follows that
$$\vp(\nu,b_0)(y)-\vp(\nu,b_0)(\sum r_ix_i)=
\sum_{\cal F}{r_{\tau,b,x}\vp(\tau,b)(x)}$$
with $r_{\tau,b,x}=0$ whenever $\tau\leq\nu.$
Hence $\vp(\nu,b_0)(y)-\vp(\nu,b_0)(\sum r_ix_i)\in \dot{M}_\nu$
as wished.\cvd\acapo

\refe
{References}
[A-B-W] K.Akin,D.A.Buchsbaum,and J.Weyman, Schur functors and
Schur complexes, {\sl Adv.in Math.}, {\bf 44} (1982), 207-278.\acapo
[B] G.Boffi, Characteristic-free decomposition of skew Schur functors,
{\sl J.Algebra}, {\bf 125} (1989), 288-297.\acapo
[C-V 1] M.Costantini,and M.Varagnolo,
 Quantum double and multiparameter quantum
groups, to appear in {\sl Comm.Algebra}.\acapo
[C-V 2] M.Costantini,and M.Varagnolo,
 Multiparameter quantum function algebra at
roots of 1, preprint 1994.\acapo
[H-H] M.Hashimoto,and T.Hayashi, Quantum
multilinear algebra, {\sl T\^ohoku Math.J.}, {\bf 44} (1992), 471-521.\acapo
[H-L-T] T.J.Hodges,T.Levasseur,and M.Toro, Algebraic structure of
multiparameter quantum \acapo groups, preprint 1994.\acapo
[N-Y-M] M.Noumi,H.Yamada,and K.Mimachi, Finite dimensional representations
of the quantum group $GL_q(n;\C)$ and the zonal spherical functions on
$U_q(n-1)\setminus U_q(n)$,preprint 1993. \acapo
[P-W] B.Parshall,and J.Wang, Quantum linear groups, {\sl Mem.of AMS},
{\bf 89} (1991), n.439.\acapo
[R] N.Y.Reshetikhin, Multiparameter quantum groups and twisted quasi-triangular
Hopf algebras, {\sl Lett.Math.Phys.}, {\bf 20} (1990), 331-335.\acapo
[S] A.Sudbery, Consistent multiparameter quantization of GL(n), {\sl
J.Phys.A}, {\bf  23} (1990), 941-961.\acapo
\end

\footnote{Note from www-admin@xxx.lanl.gov: file may have been truncated
during email transfer}

\bye